\documentclass{article}

\usepackage{arxiv}

\usepackage[utf8]{inputenc} 
\usepackage[T1]{fontenc}    
\usepackage{hyperref}       
\usepackage{url}            
\usepackage{booktabs}       
\usepackage{amsfonts}       
\usepackage{nicefrac}       
\usepackage{microtype}      
\usepackage{lipsum}
\usepackage{graphicx}
\graphicspath{ {./images/} }

\title{
Approach to Evaluating Characteristics of Multichannel Loss System with FCFD Preempted Priority Discipline
}

\author{
 Tatashev A.G. \\
  Department of Higher Mathematics\\
  Moscow Automobile and Road Construction\\
  State Technical University (MADI) \\
  Moscow, Leningradsky avenue, 64, Russia  \\
  \texttt{a-tatashev@yandex.ru} \\
   \And
Seleznjev O.V.\\
Umea University, Sweden\\
  \texttt{oleg.seleznjev@umu.se} \\
   \And
 Yashina M.V. \\
  Department of Higher Mathematics\\
  Moscow Automobile and Road Construction\\
  State Technical University (MADI) \\
  Moscow, Leningradsky avenue, 64, Russia  \\
  \texttt{mv.yashina@madi.ru} \\
}

\begin{document}
\maketitle
\begin{abstract}
In the paper, we consider a multichannel loss preemptive priority system with a Poisson input and general service time distribution depending on the priority of job. Jobs of the same priority are preempted according with First Come First Displaced (FCFD) protocol.
Approximate formulas are obtained for the loss probability of a prescribed priority job  and some other characteristics of the system.
It particular cases, the obtained formulas are exact.
\vskip 10pt

\end{abstract}

\keywords{
Queueing systems
\and
 loss multichannel system
\and
loss probability
\and
first come~--- last displaced protocol (FCLD);
\and
first come~--- first displaced protocol (FCFD)
 \and
limited processor sharing
\and  approximate formulas
}

\section{Introduction}


Exact results for characteristics of queueing systems with priority disciplines are known for one-channel system and systems in that service time
is distributed exponentially with parameter independent of job priority~\cite{Jaiswal}--~\cite{Takagi16}. Approximate approaches to evaluate characteristics of multichannel priority waiting  queueing systems with general service time distributions was proposed
 in~[8]--~\cite{Alencar20-11}.
 In ~\cite{MT-prior-92-12}, an approximate approach is proposed to evaluate characteristics of a loss multichannel preemptive priority system such that,
 if a job arrives to the system and all servers are processing jobs of the same priority or higher, then the arriving job is lost.
 If all servers are busy and at least one job with a lower priority is serviced, then the job of the lowest priority among the priorities of jobs
 in service is preempted and lost such that the preempted job arrived later than the other jobs of this priority (Last Come, First Displaced protocol, LCFD).

This paper proposes an approximate approach to evaluate the loss probability and some other characteristics of a loss multichannel preemptive priority system such that, if a job arrives to the system, and all servers are processing jobs of the higher priority, then the arriving job is lost. If all servers are busy and at least one job with the same or a lower priority is serviced, then a job of the lowest priority among the priorities of jobs in service is lost such that this job arrived earlier than the other jobs of this priority (First Come, First Displaced protocol, FCFD). The FCFD disciplines are useful when the importance of call decreases with the elapsing time~\cite{Katchner L.-13}. In the case of one-channel and in the case of exponential service time distribution, the proposed formulas are exact. The formulas are also exact if, with prescribed probabilities, the service time equals~0 or is distributed exponentially with parameter independent of the job priority. Section~2 describes the considered system. In Section~3, an approach is proposed to evaluate the probability that a job of a prescribed priority is lost at the arrival moment. Section~4 proposes approach to evaluate the probability that a accepted job of prescribed priority is lost due to preemption. An approximate formula for the loss probability is proposed in Section~5.  In Section~6, a numeric example is presented.

\section{
System Description
}
\label{section:SD}

We use the following notations: $v_i$ is the average sojourn time for a job of the priority $i;$ $w_i$ is the average waiting time including the service interruptions for the priority $i$ job; $p_i$ is the probability that the service of priority $i$ class job does not start immediately; $u_i$ is the average time before the start of the priority $i$ job service provided this time is not equal to 0; $h_i$ is the average number of the priority $i$ job preemptions $i=1,\dots,N;$ $g_i$ is the duration of a service interruption interval for the priority $i$ job, $i=2,\dots,N;$ $\Lambda_{i}$ is the total arrival rate of priority-classes no lower than $i:$ $\Lambda_i=\lambda_1+\dots+\lambda_i;$ $R_i$ is the load due to priority-classes no lower than $i:$ $R_i=(\lambda_1+\dots+\lambda_i)/m,$   $i=1,\dots,N.$

Denote by $c_i$ the probability of non-zero waiting for $M/G/m$ system computed by well-known Erlang's formula for a waiting system with the arrival load $R_i:$
$$c_i=\frac{(mR_i)^m}{m!(1-R_i)\sum\limits_{k=0}^{m-1}
\frac{(mR_i)^k}{k!}+(mR_i)^m}.$$

The  following equalities are true: $$w_1=p_1u_1,\eqno(1)$$
$$w_i=p_iu_i+h_ig_i,\ i=2,\dots,N,\eqno(2)$$
$$v_i=w_i+b_i,\  i=1,\dots,N,\eqno(3)$$
$$h_i=\frac{\Lambda_i(p_i-p_{i-1})}{\lambda_i},\ i=1,\dots,N.\eqno(4)$$
The proof of (4) is similar to the proof of an analogous statement for a preemptive priority system such that, in this system, the service distribution is exponential with average value independent of the priority class.
The proof of (4) is the following. The probability that all servers are busy by jobs of priority-classes not lower $i$ and there is at least one job of priority-class $i$ equals the difference of the probability that all servers are busy by jobs of priorities not lower than $i$ and  the probability that all servers are busy by jobs fo priorities not lower $i-1,$ and therefore this probability is $p_i-p_{i-1}.$ Hence the average number of preemptions of the priority $i$ per a time unit is equal to $\Lambda_{i-1}(p_i-p_{i-1}).$ From this, taking into account that the average number of arriving priority-class $i$ jobs per a time unit is equal to $\lambda_i,$ one gets (1).

\section{
Evaluation of probability that job is lost at its arrival
moment
}
\label{section:EP}

\hskip 18pt Suppose $\Lambda_i$ is the total arrival rate of the priority $i$ and higher, $\Lambda_i=\sum\limits_{j=1}^i \lambda_j;$ $R_i$ is the arriving load due to priority $i$ and higher; $R_i=\sum\limits_{j=1}^i \lambda_jb_j;$ $\beta(s)$ is the Laplace--Stieltjes of distribution $B_i(x),$ $\beta_i(s)=
\int\limits_0^{\infty}e^{-sx}dB_i(x);$ $q_i$ is the probability  that a job is lost at its arrival moment, $r_i$ is the probability that a job of the $j$th priority is lost provided that this job was accepted for the service; $\gamma_i$ is the loss probability for a job of $i$th priority, $j=1,\dots,i.$

Let the considered queuing system be called the system $S.$ Let us describe auxiliary $m$-channel queueing systems $S_i,$ $i=1,\dots,N.$ There are $N$ Poisson arrival processes of different priorities, and the rate of the $i$th priority processes is equal to $\lambda_j,$ $j=1,\dots,N.$ If there are less than $m$ jobs in the system, then they are serviced as usual. If there are $m-1$ jobs in the system, and a new job arrives, then the arriving job starts to be serviced at a rate increased by $m$ times, and the service of the other jobs stops until their number becomes less than $m.$ If there are $m$ jobs in the system and a new job arrives, and the priority of the serviced job is higher than the priority of the arriving job, then the arriving job is lost. If the priority of the arriving job is not lower, then the priority of the serviced job, then the arriving job starts to be serviced at a rate increased by $m$ times, and the job that was in service is lost.

Let us introduce the queuing system $S_i',$ $i=1,\dots,N.$ This system is a loss one-channel system with a Poisson input and the preemptive priority discipline. The rate of the arrival process for the $i$th priority job is equal to $\lambda_j,$ $j=1,\dots,i.$ The service time distribution for a job of the $j$th priority is
$B_j(mx),$ i.e., the service rate increases by $m$ times.

Suppose $p_{ki}$ is the stationary probability that there are $k$ jobs in the system $S_i,$
$k=0,1,\dots,m;$ $p_{ki}'$ is the stationary probability that there are $k$ jobs in the system
 $S_i',$ $k=0,1;$ $g_i$ is the average duration of the busy period of the system $S_i'.$ We have
$$p_{0i}'=\frac{1/\Lambda_i}{g_i+\frac{1}{\Lambda_i}},\  p_{1i}'=\frac{g_i}{g_i+\frac{1}{\Lambda_i}},\ i=1,\dots,N.$$

Let us obtain the value $g_i.$ Denote by $d_j$ the average duration of a busy period of the system $S_i'$ provided that this period  starts with servicing a job of the $j$th priority, $j=1,\dots,i.$ The probability that a busy period of the system $S_i'$ starts with servicing a job of the $j$th probability is $\lambda_j/\Lambda_i.$ Using the total probability formula, we get
$$g_i=\frac{1}{\Lambda_i}\sum\limits_{j=1}^i
\lambda_jd_j.\eqno(1)$$
The probability that during the service time of a job of $j$th priority, no job of a priority not lower than $i$ is
$\int\limits_0^{\infty}e^{-\Lambda_i x}dB_j(mx),$ or $\beta_j\left(\frac{\Lambda_j}{m}\right).$ If the the $j$th priority job is preemptied, then, with probability $\lambda_s/\Lambda_j,$ the preemption is due the arrival of the $s$th priority job, and  the average remaining time of the busy period equals $g_s,$ $s=1,\dots,j.$ The average service time of the $j$th priority job including service interruption time equals
$\int\limits_0^{\infty}(1-B_j(mx))e^{-\Lambda_jx}dx,$ or $\left(1-\beta_i\left(\frac{\Lambda_j}{m}\right)\right)/\Lambda_j.$ Hence,
$$d_i=\left(1-\beta_j\left(\frac{\Lambda_j}{m}\right)\right)\left(\frac{1}{\Lambda_j}+g_j\right),\ j=1,\dots,N.\eqno(2)
$$
Using (1), (2), we get
$$g_i=\frac{1}{1-\frac{\lambda_i}{\Lambda_i}             \left(1-\beta_i\left(\frac{\Lambda_i}{m}\right)\right)}\left(\frac{1}{\Lambda_i}
\sum\limits_{j=1}^{i-1}\lambda_jd_j+
\frac{\lambda_i}{\Lambda_i^2}\left(1-
\beta_i\left(\frac{\Lambda_i}{m}\right)\right)\right),\ i=1,\dots,N.\eqno(3)$$
Using recurrent procedure based on formulas (2), (3), we can compute the value $g_i.$

If the time intervals during that there $m$ jobs in the system $S_i,$ then this system is equivalent to the usual queueng system with arrival load equal to $R_i.$
Therefore, $p_{ki}=\frac{R_i^k}{k!}p_{0i},$ $k=0,1,\dots,m-1.$ Let $\nu_i(t),$ $\nu_i'(t)$ be the number
of jobs at time $t$ in the systems $S_i$ and $S_i'$ respectively. If the intervals are excluded such that, in the system $S_i,$ there are less than $m-1$ jobs, then  $\nu_i(t)-m+1$ and $\nu_i'(t)$ are equivalent stochastic processes. Therefore $\frac{p_{mi}}{p_{m-1,i}}=\frac{p_{1i}'}{p_{0i}},$ and hence,
$p_{mi}=\frac{R_i^{m-1}\Lambda_i g_ip_{i0}}
{(m-1)!}.$

From the normalizing condition $\sum\limits_{k=0}^m p_{ki}=1,$ we get
$$p_{ik}=\frac{R_i^k}{k!\left(\sum\limits_{k=0}^{m-1}\frac{R_i^k}{k!}+\frac{R_i^{m-1}\Lambda_ig_i}{(m-1)!}\right)},\ k=0,1,\dots,m-1,$$
$$p_{mi}=\frac{R_i^{m-1}\Lambda_ig_i}
{(m-1)!\sum\limits_{k=0}^{m-1}
\frac{R_i^k}{k!}+R_i^{m-1}\Lambda_ig_i}.\eqno(4)$$
 Denote by $q_i$ the probability that a job is lost at its arrival moment for the system $S.$ Assume that $p_{m,i-1}$ is the approximate value of this probability. Since the jobs arrive according to a Poisson process, the value $q_i$ is also equal to the probability that, in the system $S_{i-1},$ there are $m$ jobs of the priorities no lower than $i-1.$ Denote by $c_i$ the value $p_{mi}.$
\vskip 5pt
The equality $q_i=c_{i-1}$ is exact in the following cases.
\vskip 5pt
1. If $m=1,$ then the system $S_i$ is equivalent to the system $S$ if it is assumed that only jobs of priority not lower than $i,$ and therefore the equality $q_i=c_i$ holds.
\vskip 5pt
2. Let the service time is distributed exponentially with
parameter $\mu$ not depending on the priority of job. The presence of lower priorities in the system $S$ does not affect the servicing of higher priority jobs. Since the service time is distributed exponentially, then remaining service time is also distributed exponentially with the same parameter. Thus the probability that, in the system $S,$ there are no jobs of the priorities lower than $i$ is the same as the probability that, in related non-priority system with arriving load $R_i$ a job is lost,  i.e., according to the first Erlang formula
$$q_i=\frac{R_i^m}{m!\sum\limits_{k=0}^m
\frac{R_i^k}{k!}}.$$
In this case, $g_i=\frac{1}{m\mu}=\frac{R_i^m}{\Lambda_im},$ then (4) may be rewritten as
$$p_{mi}=\frac{R_i^m}
{m!\sum\limits_{k=0}^m
\frac{R_i^k}{k!}}.$$
Thus, $q_i=c_{i-1}.$
\vskip 5pt
3. The equality $q_i=c_{i-1}$ is also holds if, with a  probability depending on the priority of job, the service time is equal to 0, and, with the additional priority, the service time is distributed exponentially with parameter independent of the priority.

\section{
Evaluation of probability that accepted job  will be preempted
}
\label{section:EP}

\hskip 18pt Denote by $r_i$ the probability that a prescribed priority job accepted to service will be preempted.

Let us prove the following equality
$$r_i=\frac{\Lambda_i(q_i-q_{i-1})}{\lambda_i(1-q_i)}.\eqno(5)$$
An arriving job of the priority $i$ or higher preempts a job of the $i$th priority if there are $m$ jobs in the system, and at least one of these jobs is a job of the $j$th priority. Therefore, with the probability $q_i-q_{i-1},$ the system is in the state such that an arriving job of the priority not lower than $i$ preempts a job of the $i$th priority. Hence the average number of interruptions of servicing the $i$th priority jobs per a time unit equals $\Lambda(q_i-q_{i-1})r_i.$ On the other hand, this average number equals $\lambda(q_i-q_{i-1}).$ Thus we have (5).

Replacing $q_i$ with $c_i$ in (5), we get the following approximate  formula
$$r_i=\frac{\Lambda_i(c_i-c_{i-1})}{\lambda_i(1-c_i)},\ i=1,\dots,N.\eqno(6)$$

\section{
Evaluation of loss probability
}
\label{section:ELP}

\hskip 18pt We have the following equality
$$\gamma_i=q_i+(1-q_i)r_i.\eqno(7)$$
Replacing $q_i$ with $c_{i-1}$ in (7) and using (6), we get
$$\gamma_i=c_{i-1}+\frac{\Lambda_i}{\lambda_i}(c_i-c_{i-1}).\eqno(8)$$

{\it Thus we propose to compute the loss probability for a job of a prescribed priority according to formula (8), where
$$c_i=\frac{R_i^{m-1}\Lambda_ig_i}{(m-1)!\sum\limits_{k=0}^{m-1}\frac{R_i^k}{k!}+R_i^{m-1}\Lambda_ig_i},$$
$$\Lambda_i=\sum\limits_{j=1}^i\lambda_j,\ R_i=\sum\limits_{j=1}^i\lambda_j b_j,$$
$$g_i=\frac{1}{1-\frac{\lambda_i}{\Lambda_i}             \left(1-\beta_i\left(\frac{\Lambda_i}{m}\right)\right)}\left(\frac{1}{\Lambda_i}
\sum\limits_{j=1}^{i-1}\lambda_jd_j+
\frac{\lambda_i}{\Lambda_i^2}\left(1-
\beta_i\left(\frac{\Lambda_i}{m}\right)\right)\right),$$
 $$d_i=\left(1-\beta_j\left(\frac{\Lambda_j}{m}\right)\right)\left(\frac{1}{\Lambda_j}+g_j\right),\ j=1,\dots,N.
$$
}
\vskip 3pt
{\bf Remark 1.} Note that, for the related multichannel loss system with LCFD preemptied priority discipline an approximate approach is
proposed in~\cite{MT-prior-92-12}  such that, according to this approach, the loss probability $\gamma_i$ for a job of the $i$th priority is computed as
$$\gamma_i=c_i+\frac{\Lambda_{i-1}}{\lambda_i}(c_i-c_{i-1}),\ m=1,\dots,N,$$
$$c_i=\frac{R_i^{m-1}\Lambda_ig_i}{(m-1)!\sum\limits_{k=0}^{m-1}\frac{R_i^k}{k!}+R_i^{m-1}\Lambda_ig_i},$$
$$\Lambda_i=\sum\limits_{j=1}^i\lambda_j,\ R_i=\sum\limits_{j=1}^i\lambda_j b_j,\
g_i=\frac{1}{\Lambda_i}\sum\limits_{j=1}^i\lambda_jd_j,$$
$$d_j=\left(1-\beta_j\left(\frac{\Lambda_{j-1}}{m}\right)\right),\
j=2,3,\dots,N,$$
$$d_1=\frac{b_1}{m}.$$

\section{
Estimating accuracy of approach
}
\label{section:EA}

\hskip 18pt  If the job service time is exponential with parameter that may depend on priority class, then the values of the loss probability for a job of prescribed priority is
the same for the FCFD preemptive discipline and the LCFD preemptive discipline.
Assume that the job service time priority is exponential with average value depending on priority-class:
$$B_i(x)=1-e^{-\mu_i x},\  i=1,2,3,$$
$N=3,$ $m=2,$ $\mu_1 =10,$  $\mu_2 =5,$ $\mu_3=2,$
$\lambda_1=\lambda_2=\lambda_3=1.$
Then the value of the loss probabilities computed according to the proposed approach are
$$\gamma_1 = 0.0045,\  \gamma_2 = 0.060,\ q_3 = 0.36.$$
The values obtained by simulation are ~\cite{MT-prior-92-12}
$$\gamma_1 = 0.0045,\ \gamma_2 = 0.064,\  \gamma_3 = 0.32.$$

\section{
Conclusion
}
\label{section:Co}

\hskip 18pt Approximate approach is proposed to evaluate the following characteristics of a multichannel loss preemptive priority  system: the loss probability of a prescribed probability job; the loss probability of a prescribed probability job at the the arrival moment;
the probability that a job of prescribed priority is lost due to preemption. In particular cases the formulas are exact. The accuracy of the approximate approach was estimated by simulation.

\bibliographystyle{unsrt}


\end{document}